\documentclass[12pt]{article}%{report}
\input{epsf.sty}
\usepackage{epsfig}
\newcommand{\oneton}{1,\cdots,n}

\makeatletter
\let\normalequation=\equation
\def\equation{\@ifnextchar[{\subequation}{\normalequation}}
\def\subequation[#1]#2{\@ifundefined{r@#1}%
  {\def\theequation{\bf ??#2}\@warning
    {Reference `#1' on page \thepage \space
     undefined}}{\edef\@tempa{\@nameuse{r@#1}}%
    \edef\theequation{\expandafter\@car\@tempa \@nil#2}}%
  \let\@currentlabel\theequation $$}
\makeatother
\begin{document}
\begin{center}
{\large\bf On Periodic Dynamical Systems} \footnote{It is
supported
  by National Science Foundation of China 69982003 and 60074005, and also supported by Graduate
  Innovation Foundation of Fudan University.}
\end{center}

\begin{center}
 Wenlian Lu, Tianping Chen\footnote{These authors are with Lab. of
Nonlinear Mathematics Science, Institute of Mathematics, Fudan
University, Shanghai, 200433, P.R.China.}
\end{center}

{\bf Abstract}\quad In this paper, we investigate the existence
and the global stability of periodic solution for dynamical
systems with periodic interconnections, inputs and
self-inhibitions. The model is very general, the conditions are
quite weak and the results obtained are universal.

\begin{center}
\end{center}

\section{Introduction}\quad
Recurrently connected neural networks, sometimes called
Grossberg-Hopfield neural networks, are described by the following
differential equations:
\begin{eqnarray}
\frac{du_{i}(t)}{dt}=-d_{i}u_{i}(t)+\sum_{j=1}^{n}a_{ij}g_{j}(u_{j}(t))+I_{i}
\quad (i=1,\cdots,n)~~(1,1) \nonumber
\end{eqnarray}
where $g_{j}(x)$ are activation functions, $d_{i}$, $a_{ij}$ are
constants and $I_{i}$ are constant inputs.

In practice, however, the interconnections   contain asynchronous
terms in general, and the interconnection weights $a_{ij}$,
$b_{ij}$, self-inhibition $d_{i}$ and inputs $I_{i}$ should depend
on time, often periodically. Therefore, we need to discuss the
following dynamical systems with time-varying delays
\begin{eqnarray}
&&\frac{du_i}{dt}=-d_{i}(t)u_{i}(t)+\sum_{j=1}^{n}a_{ij}(t)g_j(u_j(t))\nonumber\\
&+&\sum_{j=1}^{n}b_{ij}(t)f_{j}(u_{j}(t-\tau_{ij}(t)))+I_i(t)\quad
(i=1,2,\ldots,n),~ ~(1,2) \nonumber
\end{eqnarray}
or its particular case
\begin{eqnarray}
&&\frac{du_i}{dt}=-d_{i}(t)u_{i}(t)+\sum_{j=1}^{n}a_{ij}(t)g_j(u_j(t))\nonumber\\
&+&\sum_{j=1}^{n}b_{ij}(t)f_{j}(u_{j}(t-\tau_{ij}))+I_i(t)\quad
(i=1,2,\ldots,n),~~(1.3) \nonumber
\end{eqnarray}
and the systems with distributed delays
\begin{eqnarray}
&&\frac{du_{i}(t)}{dt}=-d_{i}(t)u_{i}(t)
+\sum_{j=1}^{n}a_{ij}(t)g_{j}(u_{j}(t))\nonumber\\
&+&\sum_{j=1}^{n}b_{ij}(t)\int_{0}^{\infty}
k_{ij}(s)f_{j}(u_{j}(t-\tau_{ij}(t)-s))ds + I_{i}(t)
~(i=1,2,\ldots,n), ~~(1.4) \nonumber
\end{eqnarray}
where $d_{i}(t)>d_{i}>0$, $a_{ij}(t), b_{ij}(t), \tau_{ij}(t)>0,
I_i(t): \mathbf{R}^{+}\rightarrow \mathbf{R}$ are continuously
periodic functions with period $\omega>0$, $i,j=1,2,\ldots,n$. For
reference, see see \cite{C,G,L,M,Zheng,Zhou} and the papers cited
in these papers.

To unify models (1.2) and (1.4), we discuss the following general
model
\begin{eqnarray}
&&\frac{du_i}{dt}=-d_{i}(t)u_{i}(t)+\sum_{j=1}^{n}a_{ij}(t)g_j(u_j(t))\nonumber\\
&+&\sum_{j=1}^{n}\int_{0}^{\infty}f_{j}(u_{j}(t-\tau_{ij}(t)-s))d_{s}K_{ij}(t,s)+I_i(t)
~(i=1,2,\ldots,n), ~~(1.5) \nonumber
\end{eqnarray}
where $d_{s}K_{ij}(t,s)$, for any fixed $t\geq 0$, are
Lebesgue-Stieljies measures and satisfy
$d_{s}K_{ij}(t+\omega,s)=d_{s}K_{ij}(t,s)$, $d_{i}(t)>0$,
$a_{ij}(t), b_{ij}(t), I_i(t), \tau_{ij}(t)>0:
\mathbf{R}^{+}\rightarrow \mathbf{R}$ are continuously periodic
functions with period $\omega>0$,

The initial condition is
%$\phi(t)=[\phi_{1}(t),\cdots,\phi_{n}(t)]^{T}$, where
\begin{eqnarray}
u_{i}(s)=\phi_{i}(s) \quad  for \quad s\in(-\infty,0],
~~(1.6)\nonumber
\end{eqnarray}
where $\phi_{i}\in C(-\infty,0]$, $i=1,\cdots,n$.

It is easy to see that if $d_{s}K_{ij}(t,0)=b_{ij}(t)$ and
$d_{s}K_{ij}(t,s)=0$, for $s\ne 0$, then (1.5)  reduces to (1.2);
In addition, if $\tau_{ij}(t)=\tau_{ij}$ are constants, then
reduces to (1.3). Instead, if
$d_{s}K_{ij}(t,s)=b_{ij}(t)k_{ij}(s)ds$, then (1.5) reduces to
(1.4).

As a precondition, we assume that for system (1.5), there exists a
unique solution with every initial condition (1.6) and the
solution continuously depends on the initial data.

\section{Main Results}\quad
For the convenience, throughout this letter, we make following two
assumptions.

{\bf Assumption 1}\quad  $|g_{j}(s)|\le G_{j}|s|+C_{j}$,
$|f_{j}(s)|\le F_{j}|x|+D_{j}$, where $G_{j}>0$, $F_{j}>0$,
$C_{j}$ and $D_{j}$ are constants ($j=1,\cdots,n$).

{\bf Assumption 2}\quad  $|g_{i}(x+h)-g(x)|\le G_{i}|h|$ and
$|f_{i}(x+h)-f(x)|\le F_{i}|h|$ ($i=1,\cdots,n$).

{\bf Main Theorem}\quad Suppose that Assumption 1 is satisfied. If
there exist positive constants $\xi_{1},\xi_{2}, \cdots,\xi_{n}$
such that for all $t>0$,
\begin{eqnarray}
-\xi_{i}d_{i}(t)+\sum\limits_{j=1
}^{n}\xi_{j}G_{j}|a_{ij}(t)|+\sum\limits_{j=1}^{n}\xi_{j}F_{j}\int_{0}^{\infty}|d_{s}K
_{ij}(t,s)|<-\eta<0 \quad (i=1,2,\cdots,n).~~(2.1)\nonumber
\end{eqnarray}
Then system (1.5) has at least an $\omega-$periodic solution
$x(t)$. In addition, if Assumption 2 is satisfied and there exists
a constant $\alpha$ such that for all $t>0$,
\begin{eqnarray}
&&-\xi_{i}(d_{i}(t)-\alpha)+\sum\limits_{j=1
}^{n}\xi_{j}G_{j}|a_{ij}(t)|\nonumber\\
&+&\sum\limits_{j=1}^{n}\xi_{j}F_{j}e^{\alpha\tau_{ij}(t)}
\int_{0}^{\infty}e^{\alpha s}|d_{s}K _{ij}(t,s)|\le 0\quad
(i=1,2,\cdots,n). ~~(2.2)\nonumber
\end{eqnarray}
Then for any solution $u(t)=[u_{1}(t),\cdots,u_{n}(t)]$ of (1.5),
\begin{eqnarray}
||u(t)-x(t)||=O(e^{-\alpha t })\quad t\rightarrow\infty.
~~(2.3)\nonumber
\end{eqnarray}

{\bf Proof:}\quad Pick a constant $M$ satisfying
$M>\frac{J}{\eta}$, where
\begin{eqnarray}
J=\max_{i}\max_{t}\bigg\{\sum\limits_{j=1}^{n}|a_{ij}(t)|C_{j}
+\sum\limits_{j=1}^{n}D_{j}\int_{0}^{\infty}|d_{s}K_{ij}(t,s)|+|I_{i}(t)|\bigg\}.
~~(2.4)\nonumber
\end{eqnarray}
and let $C=C((-\infty,0],R^{n})$ be the Banach space with norm
\begin{eqnarray}
\|\phi\|=\sup\limits_{\{-\infty<\theta\le
\omega\}}\|\phi(\theta)\|_{\{\xi,\infty\}}, \nonumber
\end{eqnarray}
where
\begin{eqnarray}
\|\phi(\theta)\|_{\{\xi,\infty\}}=\max_{i=1,\cdots,n}\xi^{-1}|\phi_{i}(\theta)|,
\nonumber
\end{eqnarray}
Denote
\begin{eqnarray}
\Omega=\{x(\theta) \in C:\|x(\theta)\|\le M,
\|\dot{x}(\theta)\|\le N \}, ~~(2.5)\nonumber
\end{eqnarray}
where
\begin{eqnarray*}
N=(\alpha+\beta+\gamma)M+c
\end{eqnarray*}
and
\begin{eqnarray*}
\alpha&=&\max\limits_{i}\sup\limits_{t}|d_{i}(t)|\xi_{i}^{-1},\\
\beta&=&\max\limits_{i,j}\sup\limits_{t}|a_{ij}(t)|\xi_{i}^{-1}G_{j},\\
\gamma&=&\max\limits_{i,j}\sup\limits_{t}\int_{0}^{\infty}|d_{s}K_{ij}(t,s)|F_{j}\xi_{i}^{-1},\\
c&=&\max\limits_{i}\sup\limits_{t}|I_{i}(t)|\xi_{i}^{-1}.
\end{eqnarray*}
It is easy to check that $\Omega$ is a convex compact set.
%$\Omega\in C$ by

Now, define a map $T$ from $\Omega$ to $C$ by
\begin{eqnarray*}
\begin{array}{lll}
T:\phi(\theta)&\rightarrow x(\theta+\omega,\phi)
\end{array}
\end{eqnarray*}
where $x(t)=x(t,\phi)$ is the solution of the system (1.5) with
the initial condition $x_{i}(\theta)=\phi_{i}(\theta)$, for
$\theta\in (-\infty,0]$ and $i=\oneton$.

In the following, we will prove that $T\Omega\subset\Omega$, i.e.
if $\phi\in \Omega$, then $x\in \Omega$. To do that, we define the
following function
\begin{eqnarray}
M(t)=\sup\limits_{s\in (-\infty,0]}\|x(t+s)\|_{\{\xi,\infty\}},
~~(2.6)\nonumber
\end{eqnarray}
It is easy to see that
\begin{eqnarray}
\|x(t)\|_{\{\xi,\infty\}}\le M(t), ~~(2.7)\nonumber
\end{eqnarray}
Therefore, what we need to do is to prove $M(t)\le M$ for all
$t>0$.

Assume that $t_{0}\geq 0$ is the smallest value such that
\begin{eqnarray}
&&\|x(t_0)\|_{\{\xi,\infty\}}=M(t_0)=M, ~~(2.8)\nonumber
\end{eqnarray}
and
\begin{eqnarray}
&&\|x(t)\|_{\{\xi,\infty\}}\le M~~if~t<t_{0}, ~~(2.9)\nonumber
\end{eqnarray}
Let $i_{0}$ be an index such that
\begin{eqnarray}
\xi_{i_{0}}^{-1}|x_{i_{0}}(t)|=\|x(t_{0})\|_{\{\xi,\infty\}},
~~(2.10)\nonumber
\end{eqnarray}
Then direct calculation gives
\begin{eqnarray}
&&\bigg\{\frac{d|x_{i_{0}}(t)|}{dt}\bigg\}_{t=t_{0}} \le
sign(x_{i_{0}}(t_{0}))\bigg\{-d_{i_{0}}(t_{0})x_{i_{0}}(t_{0})
+\sum\limits_{j=1}^{n}a_{i_{0}j}(t_{0})
g_{j}(x_{j})\nonumber\\
&+&\sum\limits_{j=1}^{n}\int_{0}^{\infty}f_{j}(x_{j}(t_{0}-\tau_{i_{0}j}(t_{0})-s))
d_{s}K_{i_{0}j}(t_{0},s)+I_{i_{0}}(t_{0})\bigg\}\nonumber\\
&\le&-d_{i_{0}}|x_{i_{0}}(t_{0})| +\sum\limits_{j=1}^{n}
|a_{i_{0}j}(t)|G_{j}|x_{j}(t_{0})|\nonumber\\
&+&\sum\limits_{j=1}^{n}F_{j}\int_{0}^{\infty}|x_{j}(t_{0}-
\tau_{i_{0}j}(t_{0})-s)| |d_{s}K_{i_{0}j}(t_{0},s)|+J\nonumber\\
&\le&\bigg[-d_{i_{0}}\xi_{i_{t_{0}}} +\sum\limits_{j=1
}^{n}|a_{i_{0}j}(t_{0})|G_{j}\xi_{j}\bigg]\|x(t_{0})\|_{\{\xi,\infty\}}\nonumber\\
&+&\sum\limits_{j=1}^{n}F_{j}\xi_{j}\int_{0}^{\infty}
\|x(t_{0}-\tau_{i_{0}j}(t_{0})-s)\|_{\{\xi,\infty\}}
|d_{s}K_{i_{0}j}(t_{0},s)|+J\nonumber\\
&\le&\bigg[-d_{i_{0}}\xi_{i_{0}}+
\sum\limits_{j=1
}^{n}|a_{i_{0}j}(t_{0})|G_{j}\xi_{j}\nonumber\\
&+&\sum\limits_{j=1}^{n}F_{j}\xi_{j}\int_{0}^{\infty}
|d_{s}K_{i_{0}j}(t_{0},s)|\bigg]M(t_{0})+J\nonumber\\
&\le&-\eta M(t_{0})+J=-\eta M+J\nonumber\\
&<&0, ~~(2.11)\nonumber
\end{eqnarray}
which means that $\|x(t)\|_{\{\xi,\infty\}}$ can never exceed $M$.
Thus, $\|x(t)\|_{\{\xi,\infty\}}\le M(t)\le M$ for all $t>t_0$.
Moreover, it is easy to see that $\|\dot{x}(\theta+\omega)\|\le
N$. Therefore, $T\Omega\subset\Omega$.

%Let $\bar{\Omega}\subset \Omega$, which consists all
%equi-uniformly continuous functions in $\Omega$. Therefore, it is
%a compact convex set and $T\bar{\Omega}\subset\bar{\Omega}$.
By Brouwer fixed point theorem, there exists $\phi^{*}\in\Omega$
such that $T\phi^{*}=\phi^{*}$. Hence
$x(t,\phi^{*})=x(t,T\phi^{*})$, i.e.,
\begin{eqnarray}
x(t,\phi^{*})=x(t+\omega,\phi),~~(2.11)\nonumber
\end{eqnarray}
which is an $\omega-$periodic solution of the system (1.5).

Now, we prove that  (2.2) implies (2.3).

Let $\bar{u}(t)=[u(t)-x(t)]$, $z(t)=e^{\alpha t}\bar{u}(t)$. We
have
\begin{eqnarray}
\nonumber& &\frac{dz_i(t)}{dt}=-(d_{i}(t)-\alpha)z_{i}(t)
+e^{\alpha t}\bigg\{\sum_{j=1}^{n}a_{ij}(t)\bigg[g_j({u}_j(t))-g_j({x}_j(t))\bigg]\\
&+&\sum_{j=1}^{n}\int_{0}^{\infty}
\bigg[f_{j}(u_{j}(t-\tau_{ij}(t)-s))-
f_{j}(x_{j}(t-\tau_{ij}(t)-s))\bigg]d_{s}K_{ij}(t,s)\bigg\}.
~~(2.12)\nonumber
\end{eqnarray}
Therefore,
\begin{eqnarray}
\nonumber&&|\frac{dz_i(t)}{dt}|\le -(d_{i}(t)-\alpha)|z_{i}(t)|
+\sum_{j=1}^{n}|a_{ij}(t)|G_{j}|z_{j}(t)|\\\nonumber
&+&\sum_{j=1}^{n}F_{j}e^{\alpha\tau_{ij}(t)}\int_{0}^{\infty}e^{\alpha
s} |z_{j}(t-\tau_{ij}(t)-s))||d_{s}K_{ij}(t,s)|\\\nonumber &\le &
\bigg[-\xi_{i}(d_{i}(t)-\alpha)
+\sum_{j=1}^{n}\xi_{j}|a_{ij}(t)|G_{j}\bigg]\|z(t)\|_{\xi,\infty}\\\nonumber
&+&\sum_{j=1}^{n}\xi_{j}F_{j}e^{\alpha\tau_{ij}(t)}\int_{0}^{\infty}e^{\alpha
s} \|z_{j}(t-\tau_{ij}(t)-s))\|_{\xi,\infty}|d_{s}K_{ij}(t,s)|.
~~(2.13)\nonumber
\end{eqnarray}
By the same approach used before, we can prove that $z(t)$ is
bounded. Then $\bar{u}(t)=O(e^{-\alpha t})$. Main Theorem is
proved.

%Following corollaries are direct consequences of the Main Theorem.

In particular, let $d_{s}K_{ij}(t,0)=b_{ij}(t)$ and
$d_{s}K_{ij}(t,s)=0$, we have

{\bf Corollary 1}\quad Suppose that Assumption 1 is satisfied. If
there exist positive constants $\xi_{1},\xi_{2}, \cdots,\xi_{n}$
such that for all $t>0$,
\begin{eqnarray}
-\xi_{i}d_{i}(t)+\sum\limits_{j=1
}^{n}\xi_{j}G_{j}|a_{ij}(t)|+\sum\limits_{j=1}^{n}\xi_{j}F_{j}|b_{ij}(t)|
<0\quad (i=1,2,\cdots,n), ~~(2.14)\nonumber
\end{eqnarray}
in particular, if
\begin{eqnarray}
-\xi_{i}d_{i}+\sum\limits_{j=1
}^{n}\xi_{j}G_{j}|a_{ij}^{*}|+\sum\limits_{j=1}^{n}\xi_{j}F_{j}|b_{ij}^{*}|
<0\quad (i=1,2,\cdots,n). ~~(2.15)\nonumber
\end{eqnarray}
Then the system (1.2) or (1.3) has at least an $\omega-$periodic
solution  $x(t)$. In addition, if Assumption 2 is satisfied, and
\begin{eqnarray}
(-d_{i}(t)+\alpha)\xi_{i}+\sum\limits_{j=1}^{n}\xi_{j}G_{j}|a_{ij}(t)|
+\sum\limits_{j=1}^{n}\xi_{j}F_{j}|b_{ij}(t)|e^{\alpha\tau_{ij}}\le
0 ~(i=1,2,\cdots,n). ~~(2.16)\nonumber
\end{eqnarray}
Then for any solution $u(t)=[u_{1}(t),\cdots,u_{n}(t)]$ of (1.2)
or (1.3), we have
\begin{eqnarray}
||u(t)-x(t)||=O(e^{-\alpha t })\quad t\rightarrow\infty.
~~(2.17)\nonumber
\end{eqnarray}

Instead, if $d_{s}K_{ij}(t,s)=b_{ij}(t)k_{ij}(s)ds$, then we have

{\bf Corollary 2}\quad Suppose that Assumption 1 is satisfied. If
there exist positive constants $\xi_{1},\xi_{2}, \cdots,\xi_{n}$
such that for all $t>0$, there hold
\begin{eqnarray}
-\xi_{i}d_{i}(t)+\sum\limits_{j=1 }^{n}\xi_{j}G_{j}|a_{ij}(t)|
&+&\sum\limits_{j=1}^{n}\xi_{j}F_{j}|b_{ij}(t)|
\int_{0}^{\infty}|k_{ij}(s)|ds<-\eta<0\nonumber\\
&&(i=1,2,\cdots,n),~~(2.18)\nonumber
\end{eqnarray}
Then the system (1.4) has at least an $\omega-$periodic solution
$x(t)$.  In addition, if Assumption 2 is satisfied and
\begin{eqnarray}
-\xi_{i}(d_{i}(t)-\alpha)+\sum\limits_{j=1
}^{n}\xi_{j}G_{j}|a_{ij}(t)|&+&\sum\limits_{j=1}^{n}\xi_{j}F_{j}e^{\alpha\tau_{ij}(t)}
\int_{0}^{\infty}e^{\alpha s}|k _{ij}(t,s)|ds\le 0\nonumber\\
&&(i=1,2,\cdots,n). ~~(2.19)\nonumber
\end{eqnarray}
Then for any solution $u(t)=[u_{1}(t),\cdots,u_{n}(t)]$ of (1.4),
we have
\begin{eqnarray}
||u(t)-x(t)||=O(e^{-\alpha t })\quad t\rightarrow\infty.
~~(2.20)\nonumber
\end{eqnarray}

\section{Comparisons}\quad
In \cite{Zhou},   by using  Mawhin continuation theory, the
authors proved the folowing

{\bf Theorem A}\quad Suppose that Assumption 1 is satisfied. If
there are real constants $\epsilon>0$, $\xi_{i}>0$,
$0<\alpha_{ij}<1$, $0<\beta_{ij}<1$, $i,j=1,2\cdots,n$, such that
\begin{eqnarray}
&&\nonumber (-d_{i}+\alpha)\xi_{i}+G_{i}\bigg[\xi_{i}|a_{ii}^{*}|
+\frac{1}{2}\sum\limits_{j\ne
i}\xi_{j}|a_{ji}^{*}|^{2\alpha_{ji}}\bigg]
+\frac{1}{2}\xi_{i}\sum\limits_{j\ne
i}G_{j}|a_{ij}^{*}|^{2(1-\alpha_{ij})}\\
&+&\frac{1}{2}F_{i}\sum\limits_{j=1}^{n}
\xi_{j}|b^{*}_{ji}|^{2\beta_{ji}}e^{\alpha\tau_{ji}}+
\frac{1}{2}\xi_{i}\sum\limits_{j=1}^{n}
F_{j}|b^{*}_{ij}|^{2(1-\beta_{ij})}e^{\alpha\tau_{ij}}< 0\quad
(i=1,2,\cdots,n), ~~(3.1)\nonumber
\end{eqnarray}
where $|a^{*}_{ij}|=\sup_{\{0<t\le\omega\}}|a_{ij}(t)|<+\infty$,
$|b^{*}_{ij}|=\sup_{\{0<t\le\omega\}}|b_{ij}(t)|<+\infty$. Then
the dynamical system (1.3) has at least an $\omega$-periodic
solution $v(t)=[v_{1}(t),\cdots,v_{n}(t)]$. Instead, if Assumption
2 is satisfied. Then for any solution
$u(t)=[u_{1}(t),\cdots,u_{n}(t)]$ of (1.3),
\begin{eqnarray}
||u(t)-v(t)||=O(e^{-\alpha t })\quad t\rightarrow\infty.
~~(3.2)\nonumber
\end{eqnarray}

In paper \cite{Zheng}, the following comparison theorem was given.

{\bf Theorem B}\quad If the set of inequalities (3.1) holds. Then
there exist constants $\theta_{i}$, $i=1,\cdots,n$, such that
\begin{eqnarray}
&&
(-d_{i}+\alpha)\theta_{i}+\sum\limits_{j=1}^{n}\theta_{j}G_{j}|a_{ij}^{*}|
+\sum\limits_{j=1}^{n}\theta_{j}F_{j}|b_{ij}^{*}|e^{\alpha\tau_{ij}}<
0 ~~(3.3)\nonumber
\end{eqnarray}
But, the converse is not true.

Therefore, the conditions (3.1) are much more restrictive than
(3.3). And Theorem A is a special case of the Corollary 1.

In paper \cite{L}, the authors claimed that they investigate model
(1.2) with time-varying delays under assumption that
$\tau_{ij}(t)$ is periodic and $0\le\tau_{ij}'(t)<1$. However, if
$0<\tau_{ij}'(t)$, then $\tau_{ij}(t)$ is not periodic. Thus,
$\tau_{ij}(t)$ must be constants. The model reduces to model
(1.3). Therefore, they  investigate only model (1.3) with constant
time delays, rather than model (1.2) with time-varying delays.

Under Assumption 2 with $g_{j}(x)=f_{j}(x)$ being increasing, they
proved that if the following inequalities
\begin{eqnarray}
&& -d_{i}+\sum\limits_{j=1}^{n}G_{j}(1+d_{i}\omega)|a_{ij}^{*}|
+\sum\limits_{j=1}^{n}F_{j}(1+d_{i}\omega)|b_{ij}^{*}|< 0 \quad
(i=1,2,\cdots,n), ~~(3.4)\nonumber
\end{eqnarray}
and some other inequalities hold, then the dynamical system has at
least a periodic solution.

It is clear that this result is also a special case of Corollary
1. Moreover, their conditions are too strong.

\section{Numerical Example}
In this section, we give a numerical example to verify our Main
Theorem. Consider a delayed neural network with 3 neurons:
\begin{eqnarray*}
&&\frac{du_{1}}{dt}=-\bigg[2.51+\frac{1}{2}\sin^{2}(\pi
t)\bigg]u_{1}(t)+\sin^{2}(2\pi t)\tanh(u_{1}(t))+\cos^{2}(2\pi
t)\tanh(u_{2}(t))\nonumber\\&+&\sin^{2}(\pi t)\tanh(u_{3}(t))
+e^{-1}\sin^{2}(4\pi t)\arctan(u_{1}(t-|\sin(2\pi
t)|))\nonumber\\& +&e^{-1}\cos^{2}(4\pi
t)\arctan(u_{2}(t-\frac{\pi}{2}|\cos(2\pi
t)|))-\frac{e^{-1}}{2}\cos^{2}(\pi t)\arctan(u_{3}(t-1))+\sin(2\pi
t),\\
%\end{eqnarray*}
%\begin{eqnarray*}
&& \frac{du_{2}}{dt}=-\bigg[0.91+0.1\sin^{2}(\pi
t)+0.5\sin^{2}(4\pi t)\bigg]u_{2}(t)- 0.5\sin^{2}(2\pi
t)\tanh(u_{1}(t))\nonumber\\ &+&0.2\cos^{2}(4\pi
t)\tanh(u_{2}(t))+0.3\sin^{2}(\pi
t)\tanh(u_{3}(t))\nonumber\\
&-&0.7e^{-1}\sin^{2}(4\pi t)\arctan(u_{1}(t-|\sin(2\pi
t)|))+0.5e^{-1}\cos^{2}(2\pi
t)\arctan(u_{2}(t-\frac{\pi}{2}|\cos(2\pi
t)|))\nonumber\\
&+&0.2e^{-1}\cos^{2}(\pi t)\arctan(u_{3}(t-1))+2\cos(\pi t),\\
%\end{eqnarray*}
%\begin{eqnarray*}
&&\frac{du_{3}}{dt}=-\bigg[0.51+0.2\cos^{2}(\pi
t)+0.2\sin^{2}(2\pi t)+0.1\sin^{2}(4\pi t)\bigg]
u_{3}(t)\nonumber\\&-&0.4\cos^{2}(\pi t)\tanh(u_{1}(t))
+0.3\sin^{2}(2\pi t)\tanh(u_{2}(t))+ 0.2\cos^{2}(4\pi
t)\tanh(u_{3}(t))\nonumber\\&+&0.2e^{-1}\sin^{2}(\pi
t)\arctan(u_{1}(t-|\sin(2\pi t)|))+0.1e^{-1}\cos^{2}(2\pi
t)\arctan(u_{2}(t-\frac{\pi}{2}|\cos(2\pi
t)|))\nonumber\\&+&0.3e^{-1}\sin^{2}(4\pi
t)\arctan(u_{3}(t-1))+2\sin(2\pi t).
\end{eqnarray*}
It is easy to see that conditions (2.14) in Corollary 1 are
satisfied. However, conditions (3.1) in Theorem A used in
\cite{Zhou} and (3.4) used in \cite{L} are not satisfied. Fig 1
shows that $u_{i}(t)$, $i=1,2,3$, converges to a periodic
function, respectively.

\setlength{\unitlength}{1cm}
\begin{picture}(10,10)
\put(2, 0){\epsfig{file=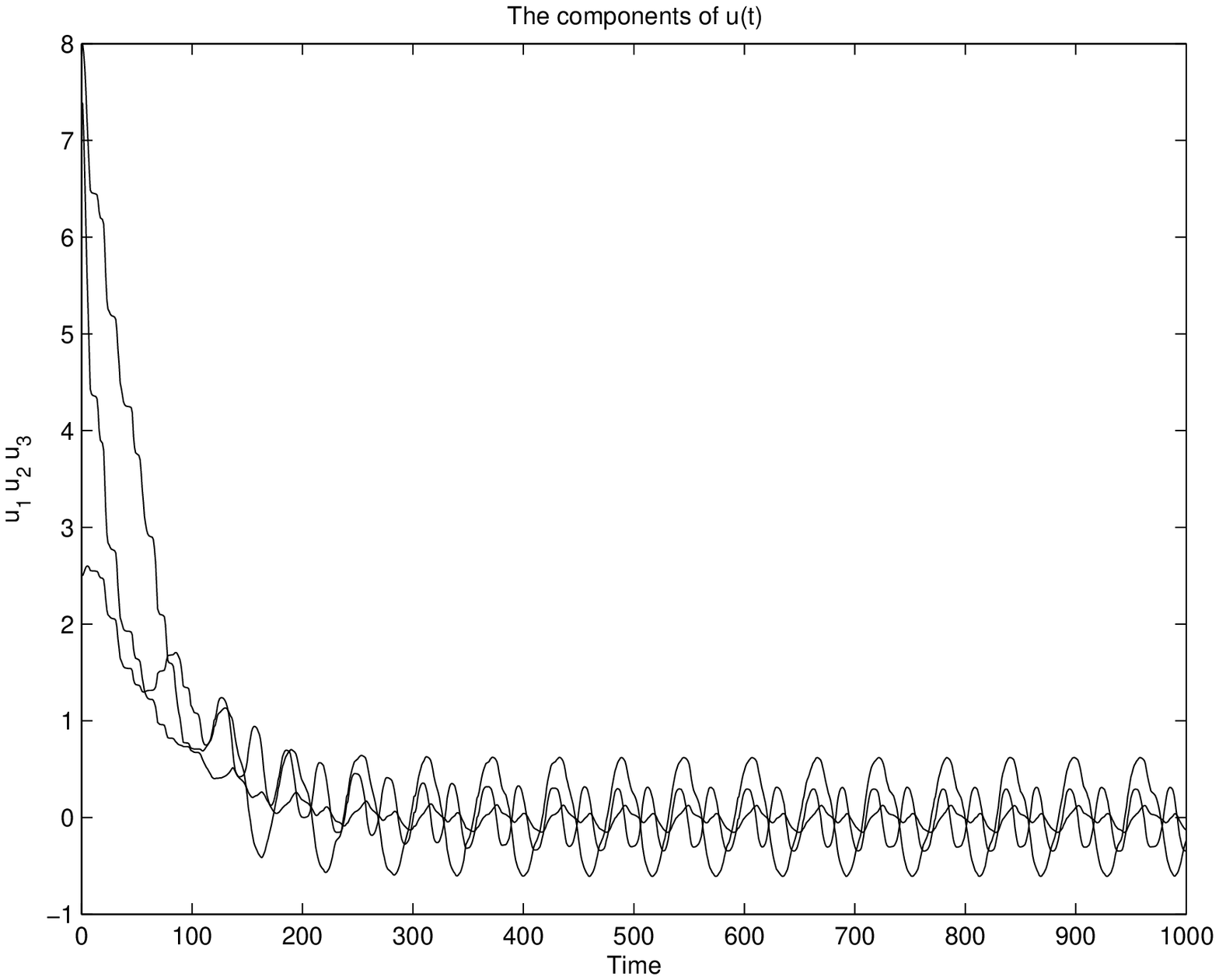,height=8cm}}
%\put(6.5,-1.0){\makebox(0,0)[c]{\footnotesize {\Large Dynamical
%Behavior of Three components} }}
\end{picture}

\section{Conclusions}\quad
In this paper, we address periodic dynamical systems. Under much
weaker conditions, the existence of periodic solution and its
exponential stability are proved.

%\setlength{\unitlength}{1cm}
%\begin{picture}(25,25)
%\put(2.5, 2.5){\epsfig{file=components.ps,height=8cm}}
%%\put(6.5,-1.0){\makebox(0,0)[c]{\footnotesize {\Large Dynamical
%%Behavior of Three components} }}
%\end{picture}

\end{document}